\documentclass[12pt]{article}

\usepackage{amstext}
\usepackage{amsthm}
\usepackage{amsmath}
\usepackage{amssymb}
\usepackage{latexsym}
\usepackage{amsfonts}
\usepackage{color}
\usepackage{graphicx}

\usepackage[pagebackref,hypertexnames=false, colorlinks, citecolor=black, linkcolor=black, urlcolor=red]{hyperref}

\usepackage[backrefs]{amsrefs}

\newcommand\black{\color{black}}

\usepackage{latexsym}
\usepackage{amssymb}
\usepackage{euscript}

\let\cal=\mathcal      

\def\mcc{M\raise.5ex\hbox{c}C}
\def\mccarthy{M\raise.5ex\hbox{c}Carthy}

\def\eg{{\it e.g. }}
\def\ie{{\it i.e. }}

\def\h{{\cal H}}

\def\K{{\cal K}}
\def\M{{\cal M}}
\def\N{{\cal N}}

\def\vare{\varepsilon}

\let\i=\infty

\def\={\ = \ }    
\def\ot{\otimes}

\def\MM{{\mathbf M}}

\def\A{{\cal A}}

\def\E{E_\l}

\def\C{\mathbb C}

\def\D{\mathbb D}

\def\dis{\displaystyle}

\def\be{\setcounter{equation}{\value{theorem}} \begin{equation}}
\def\ee{\end{equation} \addtocounter{theorem}{1}}
\def\beq{\begin{eqnarray*}}
\def\eeq{\end{eqnarray*}}
\def\se{\setcounter{equation}{\value{theorem}}} 
\def\att{\addtocounter{theorem}{1}}
\def\vs{\vskip 5pt}

\def\bp{{\sc Proof: }}
\def\ep{{}{\hfill $\Box$} \vskip 5pt \par}

\def\bl{\begin{lemma}}
\def\el{\end{lemma}}
\def\bt{\begin{theorem}}
\def\et{\end{theorem}}
\def\bprop{\begin{prop}}
\def\eprop{\end{prop}}
\def\bd{\begin{definition}}
\def\ed{\end{definition}}
\def\br{\begin{remark}}
\def\er{\end{remark}}
\def\bexer{\begin{exercise}}
\def\eexer{\end{exercise}}

\newtheorem{theorem}{Theorem}[section]
\newtheorem{prop}[theorem]{Proposition}
\newtheorem{lemma}[theorem]{Lemma}
\newtheorem{cor}[theorem]{Corollary}

\newtheorem{definition}[theorem]{Definition}
\newtheorem{examp}[theorem]{Example}

\newtheorem{defin}[theorem]{Definition}

\def\gdel{G_\delta}

\def\gdelr{G_{\delta/r}}

\def\L{{\mathcal L}}
\def\id{{\rm id}}

\def\gdels{\gdel^\sharp}

\def\vare{\varepsilon}

\def\M{{\mathbb M}}
\def\bh{{B(\h)}}
\def\bhd{\bh^d}

\renewcommand\N{{\mathbb N}}
\def\mn{\M_n}
\def\mnd{\mn^d}

\def\md{{\mathbb M}^{[d]}}

\def\d{\delta}

\def\pd{{\mathbb P}^d}

\newcommand\higd{H^\i(\gdel)}
\newcommand\bij{{\mathbb B}_{I \times J}}
\newcommand\hibij{H^\i(\bij)}
\newcommand\lkk{{\mathcal L(\K_1, \K_2)}}
\def\idcn{{\rm id}_{\mathbb{C}^n}}
\def\idX{{\rm id}_{X}}
\newcommand\dd{\d_{\MM}}

\newcommand\Fs{F^\sharp}
\newcommand\fs{f^\sharp}
\renewcommand\E{{\mathcal E}}
\renewcommand\MM{{\mathcal M}}
\newcommand\EE{\E_{\MM}}

\newcommand\cls{\rm scl_{SOT} }

\newcommand\cn{{\C^n}}
\newcommand\rd{{\tfrac{1}{r} \delta}}

\numberwithin{equation}{section}
\title{Non-commutative Functional Calculus and Spectral Theory}
\author{Jim Agler
\thanks{Partially supported by National Science Foundation Grant
DMS 1361720}
\\ U.C. San Diego\\ La Jolla, CA 92093
\and
John E. M\raise.5ex\hbox{c}Carthy
\thanks{Partially supported by National Science Foundation Grant  
DMS 1300280
}
\\ Washington University\\ St. Louis, MO 63130
}

\begin{document}

\bibliographystyle{plain}
\maketitle

\begin{abstract}
We develop a functional calculus for $d$-tuples of non-commuting elements in a Banach algebra.
The functions we apply are free analytic functions, that is nc functions that are bounded on certain
polynomial polyhedra.
\end{abstract}

\section{Introduction}

\subsection{Overview}

The purpose of this note is to develop an approach to functional calculus and spectral theory for $d$-tuples of  
elements of a Banach algebra,
 with no assumption that the elements commute.

In \cite{tay73}, J.L. Taylor considered this problem,
for $d$-tuples in $\L(X)$, the bounded linear  operators on a Banach space $X$.
 His idea was to start with the algebra $\pd$, the algebra of free polynomials\footnote{We shall use {\em free polynomial} and {\em non-commuting polynomial}
in $d$ variables interchangeably to mean an element of the algebra over the free semi-group with $d$ generators.}
 in $d$ variables over the complex numbers,
and consider what he called ``satellite algebras'', that is algebras
$\A$ that contained $\pd$, and with the property that every
representation from $\pd$ to $\L(X)$ that  extends  to a representation of $\A$
has a unique extension.
As a representation of $\pd$ is determined by choosing the images of the generators, 
\ie choosing $T = (T^1, \dots, T^d) \in
\L(X)^d$, the extension of the representation to $\A$, when it exists,
 would constitute an $\A$-functional calculus for $T$.
The class of satellite algebras that Taylor considered, which he called free analytic algebras, 
were intended to be non-commutative generalizations of the algebras $O(U)$, the algebra of holomorphic functions on a domain $U$ in
$\C^d$ (and indeed he proved in \cite[Prop 3.3]{tay73} that when $d=1$, these constitute all the
free analytic algebras). 
Taylor had already developed a successful
$O(U)$ functional calculus for $d$-tuples $T$ of {\em commuting} operators on $X$ for which a certain
spectrum (now called the Taylor spectrum) is contained in $U$ --- see \cite{tay70b, tay70a} for the original articles,
and also the article \cite{put83} by M. Putinar showing uniqueness. An excellent treatment is in \cite{cur88} by R. Curto.
However, in the non-commutative case, Taylor's approach in \cite{tay73, tay72} using homological algebra was only
partially successful.

What would constitute a successful theory? This is of course subjective, but we would argue that it should contain 
some of the following ingredients, and one  has to make trade-offs between them.
The {\em functional calculus} should  use algebras $\A$ that one knows something about --- the more the algebras are
understood, the more useful the theory. Secondly, the condition for when a given $T$ has an  $\A$-functional calculus
should be related to the way in which $T$ is presented as simply as possible. 
Thirdly,  the more explicit the map that sends $\phi$ in $\A$ to $\phi(T)$ in $\L(X)$, the easier it is to use the theory.
Finally, restricting to the commutative case, one should have a theory which agrees with the normal
idea of a functional calculus.

One does not need an explicit notion of spectrum in order to have a functional calculus.
If one does have a spectrum, it should be a collection of simpler objects than $d$-tuples in $\L(X)$,
just as in the commutative case the spectrum is a collection of $d$-tuples of complex numbers, which say something about
a commuting $d$-tuple on a Banach space.

\vs
The approach that we advocate in this note is to replace $\C^d$ as the universal set by the nc-universe
\[
\md \ :=\ \cup_{n=1}^\i \mnd ,
\]
where $\mn$ denotes the $n$-by-$n$ matrices over $\C$, with the induced operator norm
from $\ell^2_n$. 
 In other words, we look at $d$-tuples of $n$-by-$n$ matrices,
but instead of fixing $n$, we allow all values of $n$.
 We shall look at certain special open sets in $\md$.
 Let $\d$ be a matrix of free polynomials in $d$ variables, and define
\be
\label{eqa1}
\gdel \=  \{ x \in \md : \| \d(x) \| < 1 \} .
\ee
The algebras with which we shall work are algebras of the form $\higd$.
We shall define $\higd$ presently, in Definition~\ref{defab2}. For now, think of it as some sort of 
 non-commutative
analogue of the
bounded analytic functions defined on $\gdel$.
We shall develop
 conditions for a $d$-tuple in $\L(X)$
to have an $\higd$ functional calculus, in other words for a particular $T \in \L(X)^d$ to have the 
property that there is a unique extension of the polynomial functional calculus to all of 
$\higd$.



\subsection{Non-commutative functions}
\label{secab}

Let $\md = \cup_{n=1}^\i \mnd$.
A {\em graded function} defined on a subset of $\md$ is a function $\phi$ with the property that
if $x \in \mnd$, then $\phi(x) \in \mn$.
If $x \in \mnd$ and $y \in \M_m^d$, we let $x \oplus y = (x^1 \oplus y^1, \dots, x^d \oplus y^d)
\in \M_{n+m}^d$, and if $s \in \mn$ we let $sx$ (respectively $xs$) denote the tuple
$(sx^1, \dots, s x^d)$ (resp. $(x^1  s, \dots, x^d s)$).

\begin{defin}\label{defa3}
An {\em nc-function} is a graded function $\phi$ defined on a
set $D \subseteq \md$ such that

i) If $x,y, x\oplus y \in D$, then $\phi(x \oplus y) = \phi(x) \oplus \phi(y)$.

ii) If $s \in \mn$ is invertible and $x, s^{-1} x s \in D \cap \mnd$, then $\phi(s^{-1} x s) = s^{-1} \phi(x) s$.

\end{defin}

Observe that any non-commutative polynomial is an nc-function on all of $\md$.
Subject to being locally bounded with respect to an appropriate topology,
nc-functions are holomorphic \cite{hkm11b, kvv14,  amfree}, and can be thought of as bearing an analogous
relationship to non-commutative polynomials as holomorphic functions do to regular polynomials.

Nc-functions have been studied by, among others:
G.~Popescu \cite{po06,po08,po10,po11};
J. Ball, G.  Groenewald and T. Malakorn \cite{bgm06};
D. Alpay and D. Kaliuzhnyi-Verbovetzkyi \cite{akv06};
 and
J.W. Helton, I. Klep and S. McCullough \cite{hkm11a,hkm11b}
and Helton and McCullough \cite{hm12}.
We refer to the book \cite{kvv14} by  Kaliuzhnyi-Verbovetskyi  and V. Vinnikov 
on nc-functions.

We shall define matrix or operator valued nc-functions in the natural way, and use upper-case letters to denote them.
\begin{defin}
\label{defb1} Let $\K_1$ and $\K_2$ be Hilbert spaces, and $D \subseteq \md$.
We say a function $F$ is an $\lkk$-valued nc-function  on $D$ 
if 
\[
\forall_{n}\ \forall_{ x \in D \cap \mn^d}\  F(x) \in \L(\C^n \otimes \K_1,\C^n \otimes \K_2),
\]
\[
\forall_{x,y, x \oplus y \in D}\ F(x \oplus y) = F(x) \oplus F(y), \text{ and}
\]
\[
\forall_{n}\ \forall_{x \in D \cap \mn^d}\ \forall_{s \in \mn}\ s^{-1} x s \in D \implies F(s^{-1}x s) = (s^{-1}\otimes {\rm id}_{\K_1})F(x)( s \otimes {\rm id}_{\K_2}).
\]
\end{defin}

A special case of $\gdel$ in \eqref{eqa1} is when $d = I J$ and
$\delta$ is the $I$-by-$J$ rectangular matrix whose $(i,j)$ entry is the $[(i-1)J + j]^{\rm th}$ coordinate function. 
We shall give this the special name $\E$:
\[
\E( x^1, \dots, x^{IJ}) \=
\begin{pmatrix}
x^1 & x^2 & \cdots & x^J \\
x^{J+1} & x^{J+2} &\cdots & x^{2J} \\
\vdots & \vdots & \ddots & \vdots \\
x^{(I-1)J + 1} & x^{(I-1)J + 2}& \cdots & x^{IJ}
\end{pmatrix} .
\]
We shall denote the set $G_\E$ by $\bij$.
\[
\bij \=
\bigcup_{n=1}^\i \left\{ x =  (x^1, \dots, x^{IJ}) \in \mn^{IJ} \ : \
\| \E(x) \| <  1 \right\} . \]

\begin{defin}
\label{defab2}
We let $\higd$ denote the bounded nc-functions on $\gdel$, and $H_{\lkk}^\i(\gdel)$
denote the bounded $\lkk$-valued nc-functions  on $D$.
\end{defin}

 These functions were studied in 
\cite{amfree} and \cite{amfreepick}. When $\K_1 = \K_2 = \C$, we shall identify
$\higd$ with $H_{\lkk}^\i(\gdel)$.
By a matrix-valued $\higd$ function, we mean an element of some  $H_{\lkk}^\i(\gdel)$
with both $\K_1$ and $\K_2$ finite dimensional.

\section{Hilbert tensor norms}
\label{secht}

We wish to define norms on matrices of elements of $\L(X)$. If $X$ were restricted to be a Hilbert space $\h$,
there would be a natural way to do this by thinking of an $I$-by-$J$ matrix in $\L(\h)$ as a linear map from 
 the (Hilbert space) tensor product $\h \otimes \C^J$ to $\h \otimes \C^I$. We would like to do this in general.
 
 Note first that although any Banach space can be embedded in an operator space
 (see \eg \cite[Chap. 3]{pi03}), which in turn can be realized as a subset of some $\L(\h)$, we would lose the multiplicative
 structure of $\L(X)$, so that will not work in general for our purpose.

Let us recall some definitions from the theory of tensor products on Banach spaces \cite{ry02, dfs08}.
A {\em reasonable cross norm}  on the algebraic tensor product $X \ot Y$ of two Banach spaces is
a norm $\tau$ satisfying

(i) For every $x \in X,\ y \in Y$, we have
$\tau (x \ot y) = \| x \| \| y \|$.

(ii)  For every $x^* \in X^*,\ y^* \in Y^*$, we have
$\| x^*\ot y^* \|_{(X \ot Y, \tau)^*} = \| x^* \| \| y^* \|$.

A {\em uniform cross norm} is an assignment to each pair of Banach spaces $X,Y$ a reasonable cross-norm on
$X \ot Y$ such that if $ R : X_1 \to X_2$ and $S : Y_1 \to Y_2$ are bounded linear operators, then
\[
\| R \ot S \|_{ X_1 \ot Y_1 \to X_2 \ot Y_2} \ \leq \ \| R \| \| S \| .
\]
A uniform cross norm $\tau$ is {\em finitely generated} if, for every pair of Banach spaces $X, Y$
and every $u \in X \ot Y$, we have
\[
\tau (u; X\ot Y) \= \inf \{ \tau (u; M\ot N),\ u \in M \ot N,\ {\rm dim\ } M < \i, {\rm dim \ } N < \i \}.
\]
A finitely generated uniform cross norm is called a tensor norm. 
Both the injective and projective tensor products are tensor norms \cite[Prop.'s 1.2.1, 1.3.2]{dfs08}, \cite[Section 6.1]{ry02}, and there are others \cite{ry02,dfs08}.
When $\tau$ is a reasonable cross norm, we shall write $X \ot_\tau Y$ for the Banach space that is the completion of 
$X \ot Y$ with respect to the norm given by $\tau$.

\begin{defin}
\label{defht1} Let $X$ be a Banach space. A Hilbert tensor norm 
 on $X$ is an assignment of a reasonable cross norm $h$
 to $X \ot \K$ 
for every Hilbert space $\K$
 with the
 property:

  If $R : X \to X$ and $S: \K_1 \to \K_2$
are bounded linear operators,  and $\K_1$ and $\K_2$ are Hilbert spaces, then
\be
\label{eqht2}
\| R \ot S \|_{\L( X \ot_h \K_1,\ X \ot_h \K_2)} \ \leq \ \| R \|_{\L(X)} \| S \|_{\L(\K_1, \K_2)}.
\ee

%
\end{defin}
%
%
%

Any uniform cross norm is a Hilbert tensor norm, but there are others. Most importantly, if $X$ is itself a Hilbert space,
then the Hilbert space tensor product is a Hilbert tensor norm.

In what follows, we shall use $\ot$ without a subscript to denote the Hilbert space tensor product of two
Hilbert spaces,
and $\ot_h$ to denote a Hilbert tensor norm.

%
%
%
%

Let $X$ be a Banach space, and let $h$ be a Hilbert tensor norm on $X$. Let 
$R = ( R_{ij})$ be an $I$-by-$J$ matrix with entries in $\L(X)$. Then we can think of $R$
as a linear operator from $X \ot \C^J$ to $X \ot  \C^I$. We shall use $h$ to define a norm  for $R$. Formally, let $E_{ij} : \C^J \to \C^I$ be the matrix with $1$ in the $(i,j)$
slot and $0$ elsewhere. Let $\K$ be a Hilbert space.
Then we define
\se
\att
\begin{eqnarray}
\nonumber
R_{h,\K} : X \ot_h( \C^J \ot \K)  &\to& X \ot_h (\C^I \ot \K) \\
R_{h,\K} &=& \sum_{i=1}^I \sum_{j=1}^J R_{ij} \ot_h ( E_{ij} \ot {\rm id}_\K)
\label{eqht4}
\end{eqnarray}
Then we define
\be
\label{eqht5}
\| R \|_h \= \sup \{ \| R_{h,\K} \| \ : \ \K \ {\rm is\ a \ Hilbert\ space }\} ,
\ee
and (borrowing notation from the Irish use of a dot or {\em s\'eimhi\'u} for an ``h'')
\be
\label{eqht6}
\| R \|_\bullet \= \inf \{ \| R \|_h \ : \ h \ {\rm is\ a \ Hilbert\ tensor\ norm }\} .
\ee

Let us record the following lemma for future use.
\bl
\label{lemht1}
Let 
$R = ( R_{ij})$ be an $I$-by-$J$ matrix with entries in $\L(X)$.
Then
\be
\label{eqht7}
\| R \|_\bullet \ \geq \ \max_{i,j} \| R_{ij} \|_{\L(X)} .
\ee
\el
\bp
Let $B_i$ be the $1$-by-$I$ matrix with $i^{\rm th}$ entry $\id_{X}$, and the other entries the $0$ element
of $\L(X)$. Let $C_j$ be the $J$-by-$1$ column matrix, with $j^{\rm th}$ entry $\id_{X}$, and the other entries  $0$.
Let $h$ be any Hilbert tensor norm on $X$.
By \eqref{eqht2}, we have
$\| B_i \|_h $ and $\| C_j \|_h$ are $\leq 1$, and since $h$ is a reasonable cross norm we get
that they both  exactly equal $1$. We have
\beq
 \| R_{ij} \|_{\L(X)}  &\=& \| B_i R C_j \|_{\L(X)}  \\
&\leq &  \|  R  \|_{\L(X \ot_h \C^J, X \ot_h \C^I)} \\
&\leq & \| R \|_h.
\eeq
Since this holds for every $h$, we get \eqref{eqht7}.
\ep

\section{Free analytic functions}
\label{secb}

Here are some of the primary results of \cite{amfree}. When  $\delta$ is an $I$-by-$J$ rectangular matrix 
with entries in $\pd$, and $x \in \mnd$, we shall think of $\d(x)$ as an element
of $\L( \C^n \otimes \C^J , \C^n \otimes \C^I)$. If $\MM$ is a Hilbert space, we shall write
$\dd(x) $ for $\delta(x) \otimes \id_{{\MM}}$, and think of it as an element of
 $$\L( \C^n \otimes \MM^J , \C^n \otimes \MM^I) \=
\L( \C^n \otimes ( \C^J \ot \MM ) , \C^n \otimes (\C^I \ot \MM)) .
$$

\bt
\label{thmb1}
Let $\delta$ be an $I$-by-$J$ rectangular matrix of free polynomials,
and assume $\gdel$ is non-empty. Let $\K_1$ and $\K_2$ be finite dimensional Hilbert spaces.
A function $\Phi$ is in  $H_{\lkk}^\i(\gdel)$ if and only if
there is a function $F$ in  $H_{\lkk}^\i(\bij)$, with $ \| F \| \leq \| \Phi\|$,  such that
$\Phi = F \circ \delta.$
\et

\bt
\label{thmb2}
 Let $\K_1$ and $\K_2$ be finite dimensional Hilbert spaces.
 If  $F$ is in  $H_{\lkk}^\i(\bij)$ and $\| F \| \leq 1$, 
then  there exists an auxiliary Hilbert space $\MM$ and an isometry
\be
\label{eqb11}
V \=  \begin{bmatrix}A&B\\C&D\end{bmatrix} \ : \K_1 \oplus {\MM}^{(I)} \to \K_2 \oplus {\MM}^{(J)}
\ee
so that for $x \in \bij \cap \mn^d$, 
\be
\label{eqb3}
F(x) \= \idcn\otimes A +  (\idcn \otimes B) \EE(x) [ \idcn \otimes \id_{{{\MM}^{(J)}}} - (\idcn \otimes D) \EE(x) ]^{-1}
(\idcn \otimes C).
\ee
Consequently, $F$ has the  series expansion
\be
\label{eqb4}
F(x) \= \idcn \ot A +
\sum_{k=1}^\i (\idcn \ot B) \EE(x) [(\idcn \ot D) \EE(x) ]^{k-1} (\idcn \ot C),
\ee
which is  absolutely convergent on $\gdel$. 
\et
If we write $_{\cn}{A}$ for $\idcn \ot A$, then equations \eqref{eqb3} and \eqref{eqb4} 
have the more easily readable form
\se
\att
\begin{eqnarray}
\label{eqb5}
F &\= & _{\cn}{A} + {}  _{\cn}{B}\  \EE\, [ I - {}_{\cn}{D}\  \EE ]^{-1} \
{}_{\cn}{C}
\\
\att
\label{eqb6}
F(x)  &=& _{\cn}{A} +
\sum_{k=1}^{\i} {}
 _{\cn}{B}\ \EE(x) \ [_{\cn}{D} \  \EE (x) ]^{k-1}\  {}_{\cn}{C}.
\end{eqnarray}

We call \eqref{eqb3} a {\em free realization} of $F$. The isometry $V$ is  not unique,
but each term on the right-hand side of \eqref{eqb6} is a free matrix-valued polynomial, each of
whose non-zero entries is homogeneous of degree $k$.
So we can rewrite \eqref{eqb6} as
\be
\label{eqb7}
F(x) \= \sum_{k=0}^\i P_k (x) 
\ee
where each $P_k$ is a homogeneous $\L(\K_1, \K_2)$-valued free polynomial, and which satisfies
\be
\label{eqb8}
\| P_k(x) \| \ \leq \ \| x \|^k \qquad \forall x \in \bij, \forall\ k \geq 1.
\ee
These formulas (\eqref{eqb5} or \eqref{eqb7}) allow us to extend the domain of $F$ from $d$-tuples of matrices 
to $d$-tuples in $\L(X)$.
Let $X$ be a Banach space, with a Hilbert tensor norm $h$.
Let $T = (T_{ij} )$ be an $I$-by-$J$ matrix of elements of $\L(X)$. If
\be
\label{eqb12}
\| T\|_h \ < \ 1 ,
\ee
where $\| T \|_h$ is defined by \eqref{eqht5},
then we can replace $\EE(x)$  in \eqref{eqb3} by $
\sum_{i,j} T_{ij} \ot_h ( E_{ij} \ot {\rm id}_\MM) $,
and 
 get a bounded operator from
$X \ot_h \K_1$ to $X \ot_h \K_2$, provided we tensor with $\idX$.

\begin{defin}
\label{defb2}
Let $\K_1$ and $\K_2$ be finite dimensional Hilbert spaces, and 
let $F$ be a  matrix-valued nc-function on $\bij$, bounded by $1$
in norm, with a free realization given by \eqref{eqb3}, and an expansion into homogeneous  
$\L(\K_1, \K_2)$-valued free polynomials given by \eqref{eqb7}.
 Let
 $T = (T_{ij} )_{i=1,j=1}^{i=I,j=J}$ be an $I$-by-$J$ matrix of bounded operators on a Banach space $X$.
Let $h$ be a Hilbert tensor norm on $X$.
Then we define $\Fs_h(T) \in \L(X \ot_h \K_1, X \ot_h \K_2)$ to equal 
\be
\label{eqb13}
\Fs_h(T) \= \sum_{k=0}^\i P_k (T)  ,
\ee
provided that the right-hand side converges absolutely.
\end{defin}
We extend the definition of $\Fs$ to functions of norm greater than $1$ by scaling. 


The definition of $\Fs_h(T)$ may seem to depend on the choice of free realization, but in fact it does not, since
the polynomials $P_k$ do not depend on the free realization.
It does depend subtly on the choice of $h$, as $\Fs_h(T)$ is a bounded linear map in
$\L(X \ot_h \K_1, X \ot_h \K_2)$, but these are all the same if $\K_1 = \K_2 = \C$.
We shall write  $\Fs(T)$ for the dim($\K_2$)-by-dim($\K_1$) matrix
\be
\label{eqb23}
\Fs(T) \= \sum_{k=0}^\i P_k (T)  ,
\ee
which is a matrix of elements of $\L(X)$.

In the following theorem we shall write $_X A$ for $\idX\otimes_h A $, and
$T_\MM$
for  $ \sum_{i,j} T_{ij} \ot_h ( E_{ij} \ot {\rm id}_\MM) $, where we assume that $h$ is understood.

\bt
\label{thmb3} Suppose $X$ is a Banach space, and $T$ is an $I$-by-$J$ matrix
of elements of $\L(X)$. Suppose $F$ is as in Theorem \ref{thmb2}, of norm at most one.

(i) If $h$ is a Hilbert tensor norm on $X$ and  $ \| T \|_h < 1$, then 
\be
\label{eqb15}
\Fs_h(T) \=
_X A +  ({} _X B) T_\MM [ I - ( {}_X D) T_\MM ]^{-1}{}
_X C,
\ee
and
\be
\label{eqb31}
\| \Fs_h(T) \| \ \leq \ \frac{1}{1- \| T \|_h} .
\ee

(ii) If $\| T \|_\bullet < 1$, then 
\be
\label{eqb32}
\| \Fs(T) \|_\bullet  \ \leq \ \frac{1}{1- \| T \|_\bullet} .
\ee

(iii) 
 If $X$ is a Hilbert space and $\tiny{H}$ is the Hilbert space tensor product,  and
$ \| T \|_{\tiny{H}} < 1$
 then
\be
\label{eqb33}
\| \Fs_{\tiny{H}}(T) \|  \ \leq \ 1.
\ee
\et
\bp (i)
Let $\| T \|_h = r < 1$. Let us temporarily denote by $G(T)$ the right-hand side of \eqref{eqb15}.
By  \ref{eqht2},
we have $ \| {}_X D \| \leq 1$, and by \eqref{eqht5}, $\| T_\MM \| < 1$. Therefore the Neumann series 
\[
[ I - ( {}_X D) T_\MM ]^{-1} \= \sum_{k=0}^\i [  {}_X D\ T_\MM ]^{k}
\]
converges to a bounded linear operator in $\L ( X \ot_h (\C^J \ot \MM))$ of norm at most
$\tfrac{1}{1-r}$. Using \ref{eqht2} again, we conclude that 
\be
\label{eqb34}
\| G (T) \|_{\L(X \ot_h \K_1, X \ot_h \K_2)} \ \leq \
1 + \frac{r}{1-r} \= \frac{1}{1-r}.
\ee
Replacing $T$ by $e^{i\theta } T$, and integrating $G( e^{i\theta } T) $ against $e^{-ik\theta}$,
we
get, for $k \geq 1$, 
\beq
\frac{1}{2\pi}
\int_0^{2 \pi}
G( e^{i\theta } T) 
e^{-ik\theta} d \theta & \= &
 {} _X B \ T_\MM [  {}_X D \ T_\MM ]^{k-1}{}_X C \\
&=&
P_k(T),
\eeq
where $P_k$ is the homogeneous polynomial from \eqref{eqb7}.
Therefore $G(T)$ is given by the absolutely convergent series
$\sum_{k=0}^\i P_k(T)$, and hence equals $\Fs(T)$, proving \eqref{eqb15}, and, by \eqref{eqb34},
also proving \eqref{eqb31}.

(ii) This follows from the definition \eqref{eqht6}.

(iii)
Using the fact that 
$\dis \begin{bmatrix}A&B\\C&D\end{bmatrix} $ is an isometry, and equation \eqref{eqb15},
 some algebraic rearrangements give
\be
\label{eqb16}
{I - \Fs_{\tiny{H}}(T)^* \Fs_{\tiny{H}} (T) } \= 
{}_X C^* [ I - T_\MM^* \ {}_X D]^{-1} [ I - T_\MM^* T_\MM] 
  [ I - ( {}_X D) T_\MM ]^{-1}{}
_X C.
\ee
Since $\| T_\MM \| < 1$, the right-hand side of \eqref{eqb16} is positive, and so the left-hand side is also,
which means $\| \Fs_{\tiny{H}}(T) \| \leq 1$.
\ep

\vs

Suppose $\Phi(x^1, \dots, x^d)$ is in $H_{\lkk}^\i(\gdel)$.
By Theorem~\ref{thmb1}, we can write $\Phi = F \circ \d$, for some $F$ in 
$H_{\lkk}^\i(\bij)$.
Let $T = (T^1, \dots , T^d) \in \L(X)^d$.
Then $\delta(T)$ is an $I$-by-$J$ matrix with entries in $\L(X)$.
If $\| \delta(T) \|_\bullet < 1$, 
then one would like 
to define $\Phi^\sharp$
by
\be
\label{eqb10}
\Phi^\sharp (T) \= \Fs (\d (T)).
\ee
As $F$ is not unique, this raises questions about whether $\Phi^\sharp$ is well-defined.
We address this in Section~\ref{secc}.


\section{Existence of Functional Calculus}
\label{secc}

Throughout this section,  $X$ will be a Banach space, 
and  $T = (T^1, \dots, T^d)$ will be a 
$d$-tuple of bounded linear operators on $X$.

Let $\d$ be an $I$-by-$J$ matrix of free polynomials in $\pd$, and let
\[
\gdel  \=   \cup_{n=1}^\i \{ x \in \mnd : \| \d(x) \| < 1 \}.
\]
We shall say that $\gdel$ is a {\em spectral set} for $T$ if
\be
\label{eqc1}
\| p (T) \|_{\L(X)} \ \leq \ \sup_{x \in \gdel} \| p(x) \| \qquad \forall \ p \in \pd.
\ee
When $P$ is an $I$-by-$J$ matrix of polynomials, then we shall consider $P$ to be an
$\L( \C^J, \C^I)$ valued nc-function. We shall let $ \M(\pd)$ denote all (finite) matrices
of free polynomials, with the norm of $P(x)$ given as the operator norm in $\L(\C^n \otimes \C^J,
\C^n \otimes \C^I)$ where $x = (x_{ij})$ is a matrix with each $x_{ij} \in \mn$.
If 
\eqref{eqc1} holds for all matrices of polynomials, \ie
\be
\label{eqc2}
\| P (T) \|_\bullet \ \leq \ \sup_{x \in \gdel} \| P(x) \| \qquad \forall\ n,\ \forall \ P \in \M( \pd),
\ee
we shall say that $\gdel$ is a {\em complete spectral set} for $T$.
If inequalities \eqref{eqc1} or \eqref{eqc2} are true with the right-hand side multiplied by a constant $K$,
we shall say $\gdel$ is a $K$-spectral set (respectively, complete $K$-spectral set) for $T$.

\bt
\label{thmc0}
The following are equivalent.

(i) There exists $s < 1$ such that
 $G_{\d/s}$ is a $K$-spectral set for $T$.

(ii) There exists $r < 1$ such that  the map
\[
\pi: f \circ (\rd) \mapsto \fs (\rd (T))
\]
 is a well-defined bounded homomorphism from $H^\i (\gdelr)$ to
$\L(X)$ with norm less than or equal to $K$ that extends the polynomial functional calculus
on $\pd \cap H^\i (\gdelr)$.

Moreover, if these conditions hold, then
  $\pi$ is the unique extension of the evaluation homomorphism on the polynomials to a bounded homomorphism from $H^\i (\gdelr)$ to
$\L(X)$.
\et
\bp (ii) $\Rightarrow$ (i):
Let $s = r$.
Let $q \in \pd$. If $\| q \|_{\gdelr}$ is infinite, there is nothing to prove.
Otherwise,
by Theorem~\ref{thmb1}, there exists $f \in H^\i (\bij)$ such that
$q = f \circ {\rd}$ on $\gdelr$, and
\[
\| f \| \ \leq \ \| q \|_{\gdelr} .
\]
Since $\pi$ is well-defined and extends the polynomial evaluation,
\[
\pi(q) \= q(T) \= \fs (\rd (T)).
\]
Therefore
\[
\| q(T) \| \ \leq \   K \, \|f \circ {\rd}  \|_{\gdelr} \ \leq  K \, \| q \|_{\gdelr} .
\]

Now, suppose (i) holds.
 Choose $r$ in $(s,1)$.
Let $\phi \in H^\i (\gdelr)$, and assume that there are functions $f_1$ and $f_2$ in
$H^\i(\bij)$ such that
\[
\phi(x) \= f_1 \circ (\rd)(x) \= f_2 \circ (\rd) (x) \qquad \forall x \in \gdelr.
\]
Expand each $f_l$ as in \eqref{eqb7} into a series of homogeneous polynomials,
so
\[
f_l(x) \= \sum_{k=0}^\i p^l_k(x), \quad l =1,2.
\]
By \eqref{eqb8}, we have $\| p^l_k (x) \| \leq \| x\|^k .$
So
\beq
\| \sum_{k=0}^N p^1_k ( \rd (x)) -
\sum_{k=0}^N p^2_k ( \rd (x)) \|_{G_{\d/s}}
&=&
\| \sum_{k=N+1}^\i p^1_k ( \rd (x)) -
\sum_{k=N+1}^\i p^2_k ( \rd (x)) \|_{G_{\d/s}}\\
&
\leq  & 2
\sum_{k=N+1}^\i (\frac{s}{r} )^k \=
 2\ \frac{ s^{N+1}}{r^N}\frac{1}{r-s} .
\eeq
Therefore
\be
\label{eqc31}
\| \sum_{k=0}^N p^1_k ( \rd (T)) -
\sum_{k=0}^N p^2_k ( \rd (T)) \|
\ \leq \ 
2K \ \frac{ s^{N+1}}{r^N}\frac{1}{r-s}.
\ee
Therefore both series $\sum_{k=0}^\i p^1_k ( \rd (T)) $ converge to the same limit,
so $\pi(\phi)$ is well-defined.

Moreover, since $\sum_{k=0}^N p^1_k ( \rd (x)) $ converges uniformly to $\phi(x) $
on $G_{\d/s}$, we have
\[
\limsup_{N \to \i} \| \sum_{k=0}^N p^1_k \circ ( \rd ) \|_{G_{\d/s}}
\ \leq \ \| \phi \|_{G_{\d/s}} 
\ \leq \ \| \phi \|_{G_{\d/r}} .
\]
Therefore
\[
\| \pi(\phi) \|_{\L(X)} \=
\lim_{N \to \i} \| \sum_{k=0}^N p^1_k ( \rd (T)) \| 
\ \leq \ K \| \phi \|_{G_{\d/r}} .
\]
The fact that $\pi$ is a homomorphism follows from it being well defined, as if
$\phi = f \circ (\rd)$ and $\psi = g \circ (\rd)$, then $\phi \psi = (fg) \circ (\rd)$.
Finally, to show that $\pi$ extends the polynomial functional calculus, suppose $q$ is a 
free polynomial in $H^\i(\gdelr)$, so $q = f \circ (\rd)$. Expand $f(x) = \sum p_k(x)$ into its homogeneous parts.
Then $\sum_{k=0}^N p_k (\rd (x))$ converges uniformly to $q(x)$ on $G_{\d/s}$, so
since $G_{\d/s}$ is a $K$-spectral set for $T$,
\[
\pi(q)
\= 
\lim_{N \to \i} \sum_{k=0}^N p_k ( \rd (T)) \= q(T) .
\]
This last argument shows that $\pi$ is the unique continuous extension of the evaluation map on 
polynomials.
\ep

A similar result holds for complete $K$-spectral sets.

\bt
\label{thmc01}
 The following are equivalent.

(i) There exists $s < 1$ such that
 $G_{\d/s}$ is a complete $K$-spectral set for $T$.

(ii) There exists $r < 1$ such that  the map
\[
\pi: F \circ (\rd) \mapsto \Fs (\rd (T))
\]
 is a well-defined completely bounded homomorphism,
satisfying
\[
\| \Fs (\rd (T)) \|_\bullet \ \leq \ K \ \| F \circ (\rd) \|_{\gdelr}
\]
that extends the polynomial functional calculus
on $\pd \cap H^\i (\gdelr)$.

Moreover, if these conditions hold, then
  $\pi$ is the unique extension of the evaluation homomorphism on the polynomials to a bounded homomorphism from $H^\i (\gdelr)$ to
$\L(X)$.
\et

The proof is very similar to the proof of Theorem~\ref{thmc1}.
The only significant difference is that \eqref{eqc31} becomes

\[
\| \sum_{k=0}^N P^1_k ( \rd (T)) -
\sum_{k=0}^N P^2_k ( \rd (T)) \|_\bullet
\ \leq \ 
2K \ \frac{ s^{N+1}}{r^N}\frac{1}{r-s}.
\]
We apply Lemma~\ref{lemht1} to conclude that both series converge to the same limit matrix.

\begin{defin}
We shall say that $T$ has a contractive (respectively, completely contractive, bounded, completely bounded) 
$\gdel$ functional calculus if there exists $ 0 < r < 1$
such that $\gdelr$ is a spectral set (respectively, complete spectral set, $K$ spectral set, complete $K$ spectral set)
for $T$.
%
\end{defin}
%
%

\begin{remark}
Even in the case
  $d = 1$, $T \in \L(\h)$, and $\delta(x) = x$, the question of when $T$ has an $H^\i(\D)$ functional calculus
  becomes murky without the {\em a priori} requirement that $\| T \| < 1$.
By  von Neumann's inequality \cite{vonN51}, 
 $T$ will have a completely contractive $\gdel$ functional calculus  if $\| T \| < 1$.
 When $\| T \| = 1$, then $p \mapsto p(T)$ extends contractively to  $H^\i(\D)$
 if $T$ does not have a singular unitary summand \cite[Thm. III.2.3]{szn-foi2}, but to guarantee
uniqueness,  the standard extra assumption is continuity in the strong operator topology
for functions that converge boundedly almost everywhere on the unit circle \cite[Section III.2.2]{szn-foi2}.

By Rota's theorem \cite{ro60}, if $\sigma(T) \subseteq (\D)$ then
$T$ is similar to an operator which 
has a  completely contractive  $H^\i(\D)$ functional calculus.
Again, the situation becomes more delicate if $\sigma (T)$ is not required to lie in $\D$.
By Paulsen's theorem \cite{pau84},
 $T$ will have a completely bounded polynomial functional calculus if 
 and only if
$T$ is similar to a
contraction.
\end{remark}

\vs

\section{Complete spectral sets}

If  $\Phi \in H^\i(\gdel)$, and $T$ is  a $d$-tuple with $ \|\d(T) \|_\bullet < 1$,
one wants to define $\Phi^\sharp(T)$ as $F^\sharp (\d(T))$. But what if there are two different functions,
$F$ and $F_1$, both in $H^\i(\bij)$, and satisfying 
\[
\Phi(x) \= F \circ \d (x) \= F_1 \circ \d(x) \quad \forall\ x \in \gdel .
\]
How does one know that $F^\sharp (\d(T)) = F_1^\sharp (\d(T))$? 
If it doesn't, is there a ``best'' choice?

 We shall say $\gdel$ is {\em bounded} if there exists $M$ such that
\[
\| x \| \leq M , \quad \forall\ x \in \gdel .
\]
This is the same as requiring that $\pd \subseteq H^\i(\gdel)$.
A stronger condition than this is to require that the algebra generated by the $\d_{ij}$ is all of $\pd$.
\begin{defin}
We shall say that $\d$ is {\em separating} if every coordinate function $x^r, 1 \leq r \leq d$, is in the
algebra generated by the functions $\{ \d_{ij} : 1 \leq i \leq I, 1 \leq j \leq J \}$.
\end{defin}

\bt
\label{thmc1}
Assume $\|\delta (T) \|_\bullet <  1 $. 

Then 
there exists $r < 1$ such that $\gdelr$ is a complete $K$-spectral set for $T$ 
if and only if
there exists $s$ in the interval  $ (\| \d(T)\|_\bullet , 1)$ such that
whenever $F$ is a matrix-valued  $\hibij$ function, and $P$ is a matrix of free polynomials
satisfying
\be
\label{eqc3}
F \circ (\tfrac{1}{s} \delta ) (x) \= P(x) \qquad \forall\ x \in G_{(1/s)\d} ,
\ee
then 
\be
\label{eqc4}
\Fs  ( \tfrac{1}{s} \d (T) )\= P(T).
\ee


 If $\d$ is separating, then it suffices to check the  condition  for the case $P=0$.
\et
\bp 
($\Rightarrow$)
 By Theorem~\ref{thmc01}, we get \eqref{eqc3} implies \eqref{eqc4} whenever
 $\gdelr$ is a  complete K-spectral set.

($\Leftarrow$)
Suppose $\| \d(T) \|_\bullet =t  <  1$, and that $s \in ( t ,1)$ has the property that
\eqref{eqc3} implies \eqref{eqc4}. Let $r = s$;
we will show that $\gdelr$ is a complete $K$-spectral set for $T$.

Let $P$ be a matrix of polynomials; we wish to show that
\be
\label{eqc5}
\| P (T) \|_\bullet \ \leq \ K \sup \{  \| P(x) \| : x \in \mnd, \| \d(x) \| < r \} .
\ee
Without loss of generality, assume that the right-hand side of \eqref{eqc5} is finite.
By Theorem~\ref{thmb1}, we can find $F$, a matrix-valued function on $\hibij$, such that
\[
 F \circ (\tfrac{1}{r} \d) \= P \quad {\rm on\ } G_{\delta /r} 
 \]
 and
\[
\| F \|  \ \leq \ \sup \{  \| P(x) \| : x \in \mnd, \| \d(x) \|  < r \} .
\]
By \eqref{eqc4}, we have 
\[
P(T) \= \Fs   (\tfrac{1}{r} \d (T)) ,
\]
and so, by Theorem~\ref{thmb3},  \eqref{eqc5} holds, with $K = \frac{r}{r-t}$ in general.

Now, suppose that 
\be
\label{eqc6}
F \circ (\tfrac{1}{s} \delta ) \= 0 \quad {\rm on\ } G_{(1/s)\d} 
\ee
implies
\be
\label{eqc7}
\Fs  (\tfrac{1}{s} \delta  (T)) \= 0.
\ee
We wish to show that \eqref{eqc3} implies \eqref{eqc4}.
Since $\d$ is separating, there is a matrix $H$ of free polynomials such that
\[
H \circ  (\tfrac{1}{s} \delta ) (x)  \= P(x) .
\]
Then
\[
(F - H) \circ (\tfrac{1}{s} \delta ) (x)  
 \= 0 \qquad \forall\ x \in G_{\d /s} ,
\]
so by hypothesis
\[
\Fs  (\tfrac{1}{s} \delta ( (T)) \= 
H^\sharp  (\tfrac{1}{s} \delta   (T)) ,
\]
and since $H$ is a polynomial,
\[
H^\sharp  (\tfrac{1}{s} \delta  (T)) \=
H  (\tfrac{1}{s} \delta  (T)) \= P(T) ,
\]
as required.
\ep

\begin{remark} 
To just check the case $P =0$, we don't need to know that $\d$ is separating, we just
need to know that whenever a polynomial is bounded on $\gdelr$, then it is expressible as
a polynomial in the $\d_{ij}$.
%
\end{remark}

Here is a checkable condition.

\bt
\label{thmd01}
Suppose $\d(0) = 0$, and that
$T \in \L(X)^d$ has 
\[
\sup_{0 \leq r \leq 1} \| \d(r T) \|_\bullet < 1 .
\]
Then $T$ has a completely bounded $\gdel$ 
functional calculus.
\et
\bp
By Theorem~\ref{thmc1}, it is sufficient to prove that \eqref{eqc3} implies \eqref{eqc4}.
Assume \eqref{eqc3} holds, \ie
\[
F \circ (\tfrac{1}{s} \delta ) (x) 
\=
 \sum_{k=0}^\i P_k ( (\tfrac{1}{s} \d(x) )
 \= P(x) \qquad \forall\ x \in G_{(1/s)\d} .
\]
By Theorem~\ref{thmb2}, $F \circ (\tfrac{1}{s} \d) - P$ has a power series expansion in a ball
centered at $0$ in $\md$.
Since $\d(0) = 0$, 
for any $m \in \N$, the number of terms in 
 $F \circ (\tfrac{1}{s} \d)(x) - P(x)$ that are of degree $m$ in $x$ is finite. 
 
 If one expands
$
 P_k ( (\tfrac{1}{s} \d(x) )$ one gets $O( (IJ)^k)$ terms, so
if 
\[
\| \tfrac{1}{s} \d (x) \| \ < \ \frac{1}{IJ} ,
\] then
the series expansion for
\[
 \sum_{k=0}^\i P_k ( (\tfrac{1}{s} \d(x) )
 \]
 converges absolutely. We conclude therefore, by rearranging absolutely convergent series, that
 if $R$ is any $d$-tuple in $\L(X)$ satisfying $\| \d(R) \|_\bullet < \tfrac{s}{IJ} $,
 then 
 \be
 \label{eqd42}
\sum_{k=0}^\i P_k ( (\tfrac{1}{s} \d(R) ) \=  P(R).
\ee
Since $\d(0) = 0$, we can apply \eqref{eqd42} to $\zeta T$,
for all sufficiently small $\zeta$.
Now we analytically continue  to $\zeta =1$, to conclude that
\eqref{eqd42} also holds for $T$.
\ep

%
%

\section{Hilbert spaces}

 If  the $d$-tuple is in 
 $ \L(\h)^d$, it is natural to work with the Hilbert space tensor product and
 the Hilbert space norm,  instead of the norm $\| \cdot \|_\bullet$.
 Throughout this section, we will assume that $S = (S^1, \dots, S^d)$ is in $\L(\h)^d$, and all norms 
 (including those used to define spectral and complete spectral sets) will be Hilbert space norms.
Many of our earlier results go through with essentially the same proofs, but, since we can use \eqref{eqb33}
instead of \eqref{eqb32}, we get better constants.

 A sample result, proved like Theorem~\ref{thmc1}, would be:
 
 \bt
 \label{thmc2}
 Let $S \in \L(\h)^d$.
 Then there exists $r < 1$ such that
 \[
 \| F^\sharp (\tfrac{1}{r} \d(S) ) \| \ \leq \ 
 \sup \{ \| F( (\tfrac{1}{r} \d(x) )\| \ : \ x \in \gdelr \} \]
 if and only if 
 
 (i) $\| \d(S) \| < 1 $ 
 
 and
 
 (ii) whenever $F$ is a matrix-valued  $\hibij$ function with
\[
F \circ (\tfrac{1}{s} \delta ) (x) \= P(x) \qquad \forall\ x \in G_{(1/s)\d} ,
\]
then 
\[
\Fs  ( \tfrac{1}{s} \d (S) )\= P(S).
\]
 \et
 
 \vs
Example \ref{examf1} shows that condition (i) does not imply (ii) in Theorem~\ref{thmc2}.

\vs
For the remainder 
of this section, fix  an orthonormal basis  $\{ e_n \}_{n=1}^\i$ for $\h$.
Then we can naturally identify $\mn$ with the operators on $\h$ that map $\vee_{k=1}^n \{ e_k \}$
to itself, and are zero on the orthogonal complement. In this way,
$\gdel$ is a subset of $\gdels$, where
\[
\gdels \ := \ \{ S \in \L(\h)^d \ : \ \| \d(S) \| < 1 \}.
\]

Since multiplication is sequentially continuous in the strong operator topology, to get a functional calculus is is enough
to know that 
$S \in \gdels$ is the strong operator topology limit of a sequence of $d$-tuples in $\gdelr$.
For any set $A \in \bhd$, let us write
$\cls(A)$ to mean the set of tuples in $\bhd$ that are strong operator topology limits of sequences from $A$.

\bt
\label{thmd1}
Suppose 
\[
S \ \in \ \bigcup_{0 < r < 1} \cls(G_{\tfrac{1}{r} \d}) .
\]
Then $S$ has 
a completely contractive $\gdel$ functional calculus.
\et
\bp
By hypothesis, there exists a sequence $x_k \in G_{(1/t) \d}$ that converges to $S$ in the strong operator topology,
for some $t < 1$.
Therefore $\d(x_k)$ converges to $\d(S)$ S.O.T., 
so $\| \d(S) \| = r \leq t < 1$.
Let $s \in (t,1)$.
By Theorem~\ref{thmc2}, it is sufficient to prove that \eqref{eqc3}
implies \eqref{eqc4}.
As in the proof of Theorem~\ref{thmc0}, we can approximate $F$ uniformly on
$\overline{\tfrac{t}{s} \bij}$ by a sequence $Q_N$, the sum of the first $N$ homogeneous polynomials.
So for all $\vare > 0$, there exists $N_0$ such that if $N \geq N_0$ then
\se
\att
\begin{eqnarray}
\label{eqd2}
\| Q_N ( \tfrac{1}{s} \d (x_k) ) - F ( \tfrac{1}{s} \d (x_k) ) \| & \ < \ & \vare \\
\att
\label{eqd3}
{\rm and\ \ } 
\| Q_N ( \tfrac{1}{s} \d (S) ) - \Fs( \tfrac{1}{s} \d (S ) \| & \ < \ & \vare .
\end{eqnarray}
As $F ((1/s) \d) = P$ on $G_{(1/s) \d}$, inequality \eqref{eqd2} means
\be
\label{eqd4}
\forall \ N \geq N_0 \quad 
\| Q_N ( \tfrac{1}{s} \d (x_k) ) - P(x_k) \| \ <\  \vare .
\ee
Since multiplication is sequentially strong operator continuous, and $Q_N$ is a matrix of polynomials,
\be
\label{eqd5}
{\rm S.O.T.} \lim_{ k \to \i} [
Q_N ( \tfrac{1}{s} \d (x_k) ) -P(x_k)]
\=
Q_N ( \tfrac{1}{s} \d (S) ) -P(S) .
\ee
The norm of a strong operator topology sequential limit is less than or equal to the limit of the norms, so
by \eqref{eqd4}, we get from \eqref{eqd5} that
\be
\label{eqd6}
\forall \ N \geq N_0 \quad 
\| Q_N ( \tfrac{1}{s} \d (S) ) -P(S) \| \ \leq \ \vare .
\ee
Using \eqref{eqd6} in \eqref{eqd3}, 
we conclude that
\[
\|  \Fs( \tfrac{1}{s} \d (S ) -P(S) \|  \ \leq \  2 \vare,
\]
Since $\vare$ was arbitrary, we conclude that \eqref{eqc4} holds,
\ie $\Fs( \tfrac{1}{s} \d (S )  = P(S)$.
\ep

\begin{cor}
\label{cord1}
Suppose each $\d_{ij}$ is the sum of a scalar and a homogeneous polynomial of degree $1$.
Then $S$ has 
a completely contractive $\gdel$ functional calculus
if and only if $\| \d(S) \| < 1$.
\end{cor}
\bp
Let $\Pi_N$ be the projection from $\h$ onto 
 $\vee_{j=1}^n \{ e_j \}$.
 Suppose $\| \d(S) \| \leq r$.
 Let $x_N = \Pi_N S \Pi_N$.
 Then $x_N$ converges to $S$ in the strong operator topology.
 Moreover,
 \[
 \d(x_n) \=
 \Pi_N \otimes {\rm id}_{\C^I} \ \d(S) \  \Pi_N \otimes {\rm id}_{\C^J} ,
 \]
 so $\| \d(x_N) \| \leq \| \d(S) \|$.
 \ep


For Hilbert spaces, replacing completely bounded by completely contractive only changes things up to similarity.
This follows from the following theorem of V. Paulsen \cite{pau84b}:
\bt
\label{thmg1}
Let $\h$ and $\K$ be Hilbert spaces, and let $A$ be a unital subalgebra of $\L(\K)$.
Let $\rho : A \to \L(\h)$ be a completely bounded homomorphism. Then there exists 
an invertible operator $a$ on $\h$, with $\| a \| \| a ^{-1} \| = \| \rho \|_{cb} $, such that
$a^{-1} \rho( \cdot) a$  is a completely contractive homomorphism.
\et
As a consequence, we get the following.
\bt
\label{thmg2}
Let $S$ be a $d$-tuple of operators on $\h$. Then $S$ has a completely bounded $\gdel$ functional calculus if and only if there exists 
an invertible operator $a$ on $\h$ such that $R = a^{-1} S a$ has a 
completely contractive $\gdel$ functional calculus.
\et
\bp Sufficiency is clear.
For necessity, 
 suppose $0 < r < 1$, and
the map 
\[
H^\i (\gdelr) \in \Phi \ \mapsto \ \Phi(S)
\]
is a completely bounded map, with c.b. norm $K$, that extends polynomial evaluations
for polynomials that are bounded on $\gdelr$.
Then in particular, 
$\gdelr$ is a complete $K$-spectral set for $S$.
Let $\{ x_k \}_{k=1}^\i$ be a countable dense set in $\gdelr$, and
let $X = \oplus x_k$.
Then for any matrix valued function $P$, we have
\[
\| P \|_{\gdelr} \= \sup \{ \| P (x) \| : x \in \gdelr \}
\=
\| P(X) \|.
\]
By hypothesis, the map
\[
\rho : P(X) \mapsto P(S)
\]
is completely bounded, with $\| \rho \|_{cb} \leq K$.
By Paulsen's theorem \ref{thmg1}, we have there exists $a$ in $\L(\h)$ such
that 
\[
P(X) \mapsto P(a^{-1} S a)
\]
is completely contractive. Therefore $\gdelr$ is a complete spectral set for
$ a^{-1} S a$.
\ep

Remark: We don't need $\K$ to be separable, so we could have taken $X$ to be the direct sum
over all of $\gdelr$. Indeed, we could sum over all $\gdelr$ which are complete $K$ spectral sets,
and get one similarity that works for all of them.

\section{Spectrum}
\label{sech}

There are several plausible ways to define a spectrum for $T \in \L(X)^d$. 
\begin{defin}
\label{defh1}
Let
\beq
\sigma_{cc}(T) &\=& \{ x \in \md: {\rm the\ map\  } p(T) \mapsto p(x) \ {\rm is\ completely\ contractive} \}\\
\sigma_{b}(T) &\=& \{ x \in \md: {\rm the\ map\  } p(T) \mapsto p(x) \ {\rm is\ bounded} \}.
\eeq
\end{defin}

By a theorem of R. Smith \cite{sm83}; \cite[Prop 8.11]{pau02}, every bounded map from an operator algebra
into a finite dimensional algebra is completely bounded. So we conclude that
if $S \in \L(\h)^d$, then 
 $\sigma_b(S)$ is the set of all $x$
that are similar to an element of $\sigma_{cc}(S)$.

\begin{defin}
\label{defh2}
Let
\beq
\Delta_{cc}(T) &\=& \{ \d : T\  {\rm has\ a\  completely\ contractive\ } \gdel \ {\rm functional\ calculus } \}
\\
\Delta_{cb}(T) &\=& \{ \d : T\  {\rm has\ a\  completely\ bounded\ } \gdel \ {\rm functional\ calculus } \}
\\
{\rm Spec}_{cc}(T) &\=& \bigcap_{\d \in \Delta_{cc}(T)} \gdel \\
{\rm Spec}_{cb}(T) &\=&\{ x \in \mnd \ : \ \forall\ \d \in \Delta_{cb}(T), \ \exists\ y \in \gdel, {\rm \ with \ } y \ {
\rm similar\ to\ } x \}.
\eeq
\end{defin}
\begin{prop}
\label{propf1}
For every $T \in \L(X)^d$, we have
\be
\label{eqh1}
\sigma_{cc}(T) \ \subseteq \  {\rm Spec}_{cc}(T) .
\ee
The set $\sigma_{cc}(T)$ is bounded, and 
for all $n$, we have $\sigma_{cc}(T) \cap \mnd$ is compact.

If $S \in \L(\h)^d$, then 
\be
\label{eqh2}
\sigma_{b}(S)\ \subseteq \  {\rm Spec}_{cb}(S) ,
\ee
and ${\rm Spec}_{cc}(S)$ is bounded.
\end{prop}
\bp
To prove \eqref{eqh1}, 
observe that if $x \in \sigma_{cc}(T)$, then $\| \d (x) \| \leq \| \d (T) \|_\bullet$, so 
$\| \d (x) \| < 1$ whenever $T$ has a completely contractive $\gdel$ functional calculus.
Boundedness follows since if $x  \in \sigma_{cc}(T)$, then $\| x^r \| \leq \|T^r \|$ for every $r$,
and continuity shows that the set is closed at each level.

For \eqref{eqh2}, if $x \in \sigma_{b}(S)$, then by Theorem~\ref{thmg1}, $x$ is similar to some $y$ in
$ \sigma_{cc}(S)$. By \eqref{eqh1}, $x$ is  similar to an element of ${\rm Spec}_{cc}(S)$, so must lie
in ${\rm Spec}_{cb}(S)$.
To see that ${\rm Spec}_{cc}(S)$ is bounded, we can use Corollary \ref{cord1} to see that if
$x \in {\rm Spec}_{cc}(S)$, then $\| \gamma (x) \| \leq \| \gamma(S) \|$ for every
matrix $\gamma$ such that each term is of total degree at most one; in particular choosing
$\gamma(x) = ( x^r - \lambda_r)/\| S^r - \lambda^r\|$ we conclude that
\[
\| x^r - \lambda^r \| \ \leq \ \| S^r - \lambda^r \| \quad \forall\ r, \lambda^r .
\]
\ep
\begin{quest}
When is \eqref{eqh1} an equality?
\end{quest}
\vs

It is possible for $\sigma_{b}(T)$ to be empty - see Example~\ref{examf2}.

\section{Examples}
\label{secf}

\begin{examp}
\label{examd1}
\rm
Suppose 
\[
\d(x) \= ( x^1 \ \dots  \ x^d) ,
\]
a $1$-by-$d$ matrix. Then $\higd$ is the algebra of all bounded nc functions defined on the row contractions.
Functions on the row contractions were studied by Popescu in \cite{po06}.
Note that a function in $\higd$ need not have an absolutely convergent power series. 
When we expand $f \in \higd$ as in \eqref{eqb6} or \eqref{eqb7}, we get
\[
f(x) \= \sum_{k=0}^\i p_k (x) ,
\]
where each $p_k$ is a homogeneous polynomial of degree $k$, having $d^k$ terms.
Knowing merely that all the coefficients are bounded, one would need $\| x^j \| < \frac{1}{d}$
for each $j$ to conclude that the series converged absolutely.
However we do know that 
\[
\sum_{k=0}^\i \| p_k (x) \| 
\]
converges for all $x$ in $\gdel$.

By Theorem~\ref{thmd01} or Theorem~\ref{thmb3}, if $T \in \L(X)^d$ satisfies $\| \d(T) \|_\bullet < 1$, then the functional calculus 
\[
F \mapsto \Fs (T) 
\]
is a completely bounded homomorphism from $\higd$ to $\L(X)$,
with completely bounded norm at most
\[
\frac{1}{1- \| \d(T) \|_\bullet} .
\]
Any function in the multiplier algebra of the Drury-Arveson space can be extended without
increase of norm to a function in $\higd$ \cite{amfreepick}, so in particular one can then apply these 
functions to $T$.
\end{examp}

\begin{examp}
\label{examd2}
\rm
This is a similar example to \ref{examd1}.
This time, let $\d$ be the $d$-by-$d$ diagonal matrix with the coordinate functions written down
the diagonal. Then $\higd$ will be the free analytic functions defined on $d$-tuples $x$ with
$\max \| x^j \| < 1$.
Again, any function that is bounded on the commuting contractive $d$-tuples can be extended
to all of $\gdel$ without increasing its norm  \cite{amfreepick}.

Let $T = (T^1, \dots, T^d) \in \L(X)^d$. We can calculate $\| \d(T) \|_\bullet$ by observing that 
one gets a Hilbert tensor norm on
$X \otimes \ell^2_m$ if one defines
\[
\| (x_1, x_2, \dots, x_m) \| \= \sqrt{ \sum \| x_j \|^2_X} .
\]
It follows that $\| \d(T) \|_\bullet \leq \max \| T^j \|$, and since this is easily seen to be a lower bound,
we conclude
\be
\label{eqj1}
\| \d(T) \|_\bullet \= \max_{1 \leq j \leq d} \| T^j \|_{\L(X)} .
\ee
So, one gets an $\higd$ functional calculus whenever \eqref{eqj1} is less than $1$.
Let us reiterate that if $f \in \higd$ and we expand it in a power series, we have no guarantee that the
series will converge absolutely whenever the norm of each $T^j$ is less than one; we need to group the terms 
as in \eqref{eqb23}.
\end{examp}

\begin{examp}
\label{examf1}
\rm
Here is an example of a polynomial that has a different norm on $\gdel$ and $\gdels$.
Consequently, $\cls(\gdel) \neq \gdels$, and condition (i) in Theorem~\ref{thmc2} does not imply
(ii).

Let $ 0 < \vare < 0.2 $. For ease of reading, we shall use $(x,y)$ instead of $(x^1, x^2)$
to denote coordinates.
Let
\[
\d(x,y) \=
\begin{pmatrix}
\tfrac{1}{\vare} (y x - I) & 0& 0  \\
0 & \tfrac{1}{1+\vare} x &0  \\
0 &0& \tfrac{1}{1+\vare} y
\end{pmatrix} .
\]
Let $p(x) = x y - I$.

Claim: 
\se\att
\begin{eqnarray}
\label{eqf1}
\| p \|_{\gdel} & \leq & \vare + 4 \vare^2 \\
\att
\label{eqf2}
\| p \|_{\gdels} & \geq & 1.
\end{eqnarray}

\bp
Let $x \in \gdel$. Then $\| y \| < 1 + \vare$, and since
$y x$ is bounded below by $1-\vare$, we conclude that
$x$ is bounded below by $\frac{1-\vare}{1 + \vare}$. 
By this, we mean that for all vectors $v$, we have
\[
\| x v \| \ \geq \ \frac{1-\vare}{1 + \vare} \| v \| .
\]
So $x$ has an inverse $z$, and 
\[
\| z \| \ \leq \ \frac{1+\vare}{1 - \vare}.
\]
Let $e = y x - I$. Then $\| e \| < \vare$, and
\[
y \=  z  + ez.
\]
Therefore
\[
p(x) \=
x z + x e z - I \= x e z ,
\]
so
\[
\| p(x) \| \ \leq \ \frac{\vare(1+\vare)^2}{1-\vare} \ \leq \vare + 4 \vare^2 ,
\]
yielding \eqref{eqf1}.

For the second inequality, let $T = (S,S^*)$, where $S$ is the unilateral shift.
Then $\|\d(S,S^*)\| = \tfrac{1}{1+\vare} < 1$, and 
$\| p(S,S^*) \| = 1$, yielding \eqref{eqf2}.
\ep
\end{examp}

\begin{examp}
\label{examf2}
\rm
It is easy for
 $\sigma_{b}(T)$ to be empty.
For example, suppose $q(x) = x^1 x^2 - x^2 x^1 - I$, and choose $T \in \L(\h)^2$ so that 
$\| q(T) \| = \tfrac{1}{2}$. (This can be done, since by \cite{bp65} any operator in $\L(\h)$ that is not a non-zero scalar plus a compact is a commutator).
Then for any $x \in \M^{[2]}$, we have
$\| q(x) \| \geq 1$, so $x \notin \sigma_{cc}(T)$.
Consequently, $\sigma_{b}(T)$ is also empty.
Note that in this example, $G_q$ is empty, though $T \in G_q^\sharp$.
\end{examp}

\begin{examp}
\label{examf3}
\rm
This is an example of our non-commutative approach applied to a single matrix.
Let 
\[
U \= \{ z \in \C : | z| < 1,\ {\rm and\ }|z-1| < 1 \}.
\]
Let $X$ be a finite dimensional Banach space, and $T \in \L(X)$ have $\sigma(T) \subset U$.
Let 
\[
\d(x) \=
\begin{pmatrix}
x & 0 \\
0 & x-1 
\end{pmatrix} .
\]
Then $\higd$ will be a space of analytic functions on $U$, but the norm will not be the sup-norm;
it will be the larger norm given by
\[
\| \phi \| \ := \ \sup \{ \| \phi(S) \| : S \in \L(\h), \| \d(S) \| < 1 \}.
\]
Indeed, by Theorem~\ref{thmb1}, the norm can obtained as
\[
\| \phi \| = \inf \{ \| g \|_{H^\i(\D^2)} \ |\  g( z, z-1) = \phi (z) \ \forall\ z \ \in U \}.
\]
(It is sufficient to calculate the norm of $g$ in the commutative case, since it always
has an extension of the same norm to the non-commutative space, by \cite{amfreepick}).

By \cite[Thm. 4.9]{agmc14a},  every function analytic on a neighborhood of $\overline{U}$ is in $\higd$.
Since $X$ is finite dimensional, $T$ is similar to an operator on a Hilbert space, and by the results of Smith and Paulsen,
this can be taken to have $U$ as a complete spectral set.

Putting all this together, we can write $T$ as $a^{-1} S a$, where $S$ is a Hilbert space operator with
$\| \d(S) \| < 1$. For any $\phi$ in $H^\i(\gdel)$, we find a $g$ of minimal norm in $H^\i(\D^2)$
such that
\[
 g( z, z-1) = \phi (z)
  \ \forall\ z \ \in U .
  \]
Finally, we get the estimate
\[
\| \phi(T) \|_{\L(X)}
\ \leq \ \| a^{-1} \| \| a \|  \| g \|_{H^\i(\D^2)}  .
\]
If we know $\max ( \| T \|, \| T -1 \|) = r < 1$, we have the estimate (which works
even if $X$ is infinite dimensional)
\[
\| \phi(T) \|_{\L(X)}
\ \leq \ \frac{1}{1-r} \| g \|_{H^\i(\D^2)}  .
\]
\end{examp}

\black

\bibliography{../../references}

\begin{bibdiv}
\begin{biblist}

\bib{agmc14a}{article}{
      author={Agler, J.},
      author={M\raise.45ex\hbox{c}Carthy, J.E.},
       title={Operator theory and the {Oka} extension theorem},
     journal={Hiroshima Math. J.},
        note={To appear},
}

\bib{amfreepick}{unpublished}{
      author={Agler, J.},
      author={M\raise.45ex\hbox{c}Carthy, J.E.},
       title={Pick interpolation for free holomorphic functions},
        note={Amer. J. Math, to appear},
}

\bib{amfree}{article}{
      author={Agler, J.},
      author={M\raise.45ex\hbox{c}Carthy, J.E.},
       title={Global holomorphic functions in several non-commuting variables},
        date={2015},
     journal={Canad. J. Math.},
      volume={67},
      number={2},
       pages={241\ndash 285},
}

\bib{akv06}{incollection}{
      author={Alpay, D.},
      author={Kalyuzhnyi-Verbovetzkii, D.~S.},
       title={Matrix-{$J$}-unitary non-commutative rational formal power
  series},
        date={2006},
   booktitle={The state space method generalizations and applications},
      series={Oper. Theory Adv. Appl.},
      volume={161},
   publisher={Birkh\"auser},
     address={Basel},
       pages={49\ndash 113},
}

\bib{bgm06}{incollection}{
      author={Ball, Joseph~A.},
      author={Groenewald, Gilbert},
      author={Malakorn, Tanit},
       title={Conservative structured noncommutative multidimensional linear
  systems},
        date={2006},
   booktitle={The state space method generalizations and applications},
      series={Oper. Theory Adv. Appl.},
      volume={161},
   publisher={Birkh\"auser},
     address={Basel},
       pages={179\ndash 223},
}

\bib{bp65}{article}{
      author={Brown, Arlen},
      author={Pearcy, Carl},
       title={Structure of commutators of operators},
        date={1965},
        ISSN={0003-486X},
     journal={Ann. of Math. (2)},
      volume={82},
       pages={112\ndash 127},
      review={\MR{0178354 (31 \#2612)}},
}

\bib{cur88}{incollection}{
      author={Curto, R.E.},
       title={Applications of several complex variables to multiparameter
  spectral theory},
        date={1988},
   booktitle={Surveys of some recent results in operator theory},
      editor={Conway, J.B.},
      editor={Morrel, B.B.},
   publisher={Longman Scientific},
     address={Harlow},
       pages={25\ndash 90},
}

\bib{dfs08}{book}{
      author={Diestel, Joe},
      author={Fourie, Jan~H.},
      author={Swart, Johan},
       title={The metric theory of tensor products},
   publisher={American Mathematical Society, Providence, RI},
        date={2008},
        ISBN={978-0-8218-4440-3},
        note={Grothendieck's r{\'e}sum{\'e} revisited},
}

\bib{hkm11a}{article}{
      author={Helton, J.~William},
      author={Klep, Igor},
      author={McCullough, Scott},
       title={Analytic mappings between noncommutative pencil balls},
        date={2011},
     journal={J. Math. Anal. Appl.},
      volume={376},
      number={2},
       pages={407\ndash 428},
}

\bib{hkm11b}{article}{
      author={Helton, J.~William},
      author={Klep, Igor},
      author={McCullough, Scott},
       title={Proper analytic free maps},
        date={2011},
     journal={J. Funct. Anal.},
      volume={260},
      number={5},
       pages={1476\ndash 1490},
}

\bib{hm12}{article}{
      author={Helton, J.~William},
      author={McCullough, Scott},
       title={Every convex free basic semi-algebraic set has an {LMI}
  representation},
        date={2012},
     journal={Ann. of Math. (2)},
      volume={176},
      number={2},
       pages={979\ndash 1013},
}

\bib{kvv14}{book}{
      author={Kaliuzhnyi-Verbovetskyi, Dmitry~S.},
      author={Vinnikov, Victor},
       title={Foundations of free non-commutative function theory},
   publisher={AMS},
     address={Providence},
        date={2014},
}

\bib{pau84b}{article}{
      author={Paulsen, Vern~I.},
       title={Completely bounded homomorphisms of operator algebras},
        date={1984},
        ISSN={0002-9939},
     journal={Proc. Amer. Math. Soc.},
      volume={92},
      number={2},
       pages={225\ndash 228},
         url={http://dx.doi.org/10.2307/2045190},
      review={\MR{754708 (85m:47049)}},
}

\bib{pau84}{article}{
      author={Paulsen, V.I.},
       title={Every completely polynomially bounded operator is similar to a
  contraction},
        date={1984},
     journal={J. Funct. Anal.},
      volume={55},
       pages={1\ndash 17},
}

\bib{pau02}{book}{
      author={Paulsen, V.I.},
       title={Completely bounded maps and operator algebras},
   publisher={Cambridge University Press},
     address={Cambridge},
        date={2002},
}

\bib{pi03}{book}{
      author={Pisier, Gilles},
       title={Introduction to operator space theory},
      series={London Mathematical Society Lecture Note Series},
   publisher={Cambridge University Press, Cambridge},
        date={2003},
      volume={294},
        ISBN={0-521-81165-1},
         url={http://dx.doi.org/10.1017/CBO9781107360235},
}

\bib{po06}{article}{
      author={Popescu, Gelu},
       title={Free holomorphic functions on the unit ball of {$B({\cal
  H})^n$}},
        date={2006},
     journal={J. Funct. Anal.},
      volume={241},
      number={1},
       pages={268\ndash 333},
}

\bib{po08}{article}{
      author={Popescu, Gelu},
       title={Free holomorphic functions and interpolation},
        date={2008},
     journal={Math. Ann.},
      volume={342},
      number={1},
       pages={1\ndash 30},
}

\bib{po10}{article}{
      author={Popescu, Gelu},
       title={Free holomorphic automorphisms of the unit ball of {$B({\cal
  H})^n$}},
        date={2010},
     journal={J. Reine Angew. Math.},
      volume={638},
       pages={119\ndash 168},
}

\bib{po11}{article}{
      author={Popescu, Gelu},
       title={Free biholomorphic classification of noncommutative domains},
        date={2011},
     journal={Int. Math. Res. Not. IMRN},
      number={4},
       pages={784\ndash 850},
}

\bib{put83}{article}{
      author={Putinar, M.},
       title={Uniqueness of {Taylor's} functional calculus},
        date={1983},
     journal={Proc. Amer. Math. Soc.},
      volume={89},
       pages={647\ndash 650},
}

\bib{ro60}{article}{
      author={Rota, Gian-Carlo},
       title={On models for linear operators},
        date={1960},
     journal={Comm. Pure Appl. Math.},
      volume={13},
       pages={469\ndash 472},
}

\bib{ry02}{book}{
      author={Ryan, Raymond~A.},
       title={Introduction to tensor products of {B}anach spaces},
      series={Springer Monographs in Mathematics},
   publisher={Springer-Verlag London, Ltd., London},
        date={2002},
        ISBN={1-85233-437-1},
         url={http://dx.doi.org/10.1007/978-1-4471-3903-4},
}

\bib{sm83}{article}{
      author={Smith, R.R.},
       title={Completely bounded maps between {C*-algebras}},
        date={1983},
     journal={J. Lond. Math. Soc.},
      volume={27},
       pages={157\ndash 166},
}

\bib{szn-foi2}{book}{
      author={Szokefalvi-Nagy, B.},
      author={Foia\c{s}, C.},
       title={Harmonic analysis of operators on {Hilbert} space, second
  edition},
   publisher={Springer},
     address={New York},
        date={2010},
}

\bib{tay70b}{article}{
      author={Taylor, J.L.},
       title={The analytic functional calculus for several commuting
  operators},
        date={1970},
     journal={Acta Math.},
      volume={125},
       pages={1\ndash 38},
}

\bib{tay70a}{article}{
      author={Taylor, J.L.},
       title={A joint spectrum for several commuting operators},
        date={1970},
     journal={J. Funct. Anal.},
      volume={6},
       pages={172\ndash 191},
}

\bib{tay72}{article}{
      author={Taylor, J.L.},
       title={A general framework for a multi-operator functional calculus},
        date={1972},
        ISSN={0001-8708},
     journal={Advances in Math.},
      volume={9},
       pages={183\ndash 252},
      review={\MR{0328625 (48 \#6967)}},
}

\bib{tay73}{article}{
      author={Taylor, J.L.},
       title={Functions of several non-commuting variables},
        date={1973},
     journal={Bull. Amer. Math. Soc.},
      volume={79},
       pages={1\ndash 34},
}

\bib{vonN51}{article}{
      author={von Neumann, J.},
       title={{Eine Spektraltheorie f\"ur allgemeine Operatoren eines
  unit\"aren Raumes}},
        date={1951},
     journal={Math. Nachr.},
      volume={4},
       pages={258\ndash 281},
}

\end{biblist}
\end{bibdiv}

\end{document}